\documentclass[11pt]{amsart}
\usepackage{amssymb}
\numberwithin{equation}{section}

\newtheorem{theorem}{Theorem}[section]

\newtheorem{corollary}[theorem]{Corollary}

\newtheorem{lemma}[theorem]{Lemma}
\theoremstyle{definition}
\newtheorem{definition}[theorem]{Definition}
\theoremstyle{remark}
\newtheorem{remark}[theorem]{Remark}

\newcommand{\A}{\mathcal{A}}

\newcommand{\R}{\mathbb{R}}

\newcommand\lie[1]{\mathfrak{#1}}

\newcommand{\fg}{\lie{g}}

\def    \inv    {^{-1}}

\newcommand{\SP}        {\operatorname{Sp}}

\newcommand{\bPsi} {\bar{\Psi}}
\newcommand{\bM}    {\bar{M}}
\newcommand{\bgamma}    {\bar{\gamma}}
\newcommand{\bsigma}    {\bar{\sigma}}
\newcommand{\btau}    {\bar{\tau}}

\begin{document}

\title{A convexity theorem for torus actions on contact manifolds}
\author{Eugene Lerman}
\address{Department of
Mathematics, University of Illinois, Urbana, IL 61801}
\email{lerman@math.uiuc.edu}
%\date{12/03/2000}
\thanks{Partially supported by the NSF grant DMS - 980305 and 
the American Institute of Mathematics.}

\begin{abstract}
We show that the image cone of  a moment map for an action of a
torus on a contact compact connected manifold is a convex polyhedral
cone and that the moment map has connected fibers provided the
dimension of the torus is bigger than 2 and that no orbit is tangent
to the contact distribution.  This may be considered as a version of
the Atiyah - Guillemin - Sternberg convexity theorem for torus actions
on symplectic cones and as a direct generalization of the convexity
theorem of Banyaga and Molino for completely integrable torus actions
on contact manifolds.
\end{abstract}

\maketitle

\section{Introduction}

The goal of the paper is to prove a convexity theorem for torus
actions on contact manifolds. Recall that a {\bf contact form} on a
manifold $M$ of dimension $2n+1$ is a 1-form $\alpha$ such that
$\alpha \wedge d\alpha^n \not = 0$.  A (co-oriented) {\bf contact
structure} on a manifold $M$ is a subbundle $\xi $ of the tangent
bundle $TM$ which is given as the kernel of a contact form.  Note that
if $f$ is any nowhere vanishing function and $\alpha$ is a contact
form, then $\ker \alpha = \ker f\alpha$.  Thus a co-oriented contact
structure is a conformal class of contact forms.  One can show that a
hyperplane subbundle $\xi$ of $TM$ is a co-oriented contact structure
if and only if its annihilator $\xi^\circ$ in $T^*M$ is a trivial line
bundle and $\xi^\circ \smallsetminus 0$ is a symplectic submanifold of
the punctured cotangent bundle $T^*M \smallsetminus 0$ (we use $0$ as
a shorthand for the image of the zero section).  In fact the map $\psi
_\alpha: M\times \R \to \xi ^\circ$, $(m,t) \mapsto t\alpha _m$
defines a trivialization, and the pull-back by $\psi_\alpha$ of the
tautological 1-form on $T^*M$ is $d(t\alpha)$.  The symplectic
manifold $(M\times (0, \infty), d(t\alpha))$ is called the {\bf
symplectization} of $(M, \alpha)$.  

Recall that a {\bf symplectic cone} is a symplectic manifold $(N,
\omega)$ with a proper action of the real line which expands the
symplectic form exponentially. For example, the action of $\R$ on
$M\times (0, \infty)$ given by $s\cdot (m, t) = (m, e^s t)$ makes the
symplectization $(M\times (0, \infty), d(t\alpha))$ of $(M, \alpha)$
into a symplectic cone.  Conversely a symplectic cone is the
symplectization of a contact manifold.

Throughout the paper $\alpha $ will always denote a contact form and
$\xi$ will always denote a co-oriented contact structure.  We will
refer either to a pair $(M, \alpha)$ or to a pair $(M, \xi)$ as a
contact manifold.

 An action of a Lie group $G$ on a contact manifold
$(M, \xi)$ is {\bf contact} if the action preserves the contact
structure.  It is not hard to show that if additionally the action of
$G$ is proper (for example if $G$ is compact) and preserves the
co-orientation of $\xi$ (for example if $G$ is connected), then it
preserves a contact form $\alpha$ with $\xi = \ker \alpha$ (see
\cite{IJL}).

\subsection*{Contact moment maps}
We now recall the notion of a moment map for an action of a group on a
contact manifold.  An action of a Lie group $G$ on a manifold $M$
naturally lifts to a Hamiltonian action on the cotangent bundle
$T^*M$.  The corresponding moment map $\Phi :T^*M \to \fg^*$ is given
by 
\begin{equation}\label{mmeq}
\langle \Phi (q, p), A\rangle = \langle p, A_M (q)\rangle,
\end{equation}
for all vectors $A\in \fg$, all points $q\in M$ and all covectors
$p\in T^*_q M$.  Here and elsewhere in the paper $A_M$ denotes the
vector field induced on $M$ by $A\in \fg$.  

If the action of the Lie group $G$ on the manifold $M$ preserves a
contact distribution $\xi$, then the lifted action preserves the
annihilator $\xi^\circ \subset T^*M$.  Moreover, if the action of $G$
preserves a co-orientation of $\xi$ then it preserves the two
components of $\xi^\circ \smallsetminus 0$.  Denote one of the components by
$\xi^\circ _+$.  In this case we {\bf define} the {\bf moment map}
$\Psi$ for the action of $G$ on $(M, \xi)$ to be the restriction of
$\Phi$ to $\xi^\circ _+$:
$$
\Psi = \Phi |_{\xi^\circ _+}.
$$ 
An invariant contact form $\alpha$ on $M$ defining the contact
distribution $\xi$ is a nowhere zero section of $\xi^\circ \to M$.  We
may assume that $\alpha (M) \subset \xi^\circ _+$.  In this case we
get a map $\Psi _\alpha :M\to \fg^*$ by composing $\Psi$ with $\alpha$:
$\Psi_\alpha = \Psi \circ \alpha$.  It follows from (\ref{mmeq}) that
\begin{equation} \label{eq1}
\langle \Psi_\alpha (x) , A\rangle = \alpha _x (A_M (x))
\end{equation}
 for all $x\in M$ and all $A\in \fg$.  Recall that the choice of a
contact form on $M$ establishes a bijection between the space of smooth
functions on $M$ and the space of contact vector fields.  It is easy
to check that for any $A\in \fg$ the contact vector field
corresponding to the function $\langle \Psi_\alpha, A \rangle$ is
$A_M$.  Thus it makes sense to think of $\Psi _\alpha$ as the moment
map defined by the contact form $\alpha$ and of $\Psi$ as the moment
map defined by the contact distribution $\xi$.  Similarly the image
$\Psi_\alpha (M)$ depends on the action and the contact form while the
image $\Psi (\xi^\circ _+)$ depends only on the action and the contact
distribution.  Clearly the two sets are related: $$
\Psi (\xi^\circ _+  ) = \R^+ \Psi_\alpha (M).
$$
\begin{definition}
Let $(M, \xi)$ be a co-oriented contact manifold with an action of a
Lie group $G$ preserving the contact structure $\xi$ and its
co-orientation.  Let $\xi^\circ_+$ denote a component of $\xi^\circ
\smallsetminus 0$, the annihilator of $\xi$ minus the zero section.
Let $\Psi: \xi^\circ_+ \to \fg^*$ denote the corresponding moment map.
The {\bf moment cone} $C(\Psi)$ is the set $$ C(\Psi) : = \Psi
(\xi^\circ_+)\cup \{0\}.  $$ Note that if $\alpha $ is an invariant
contact form with $\xi = \ker
\alpha$ and $\alpha (M) \subset \xi^\circ_+$,  and if $\Psi_\alpha :M\to
\fg^*$ is the moment map defined by $\alpha$ then $C(\Psi ) = \{t f
\mid f \in \Psi _\alpha (M), \, t \in [0, \infty)\}$.
\end{definition}

We can now state the main result of the paper.

\begin{theorem}\label{main theorem}
Let $(M, \xi)$ be a co-oriented contact manifold with an effective
action of a torus $G$ preserving the contact structure and its
co-orientation. Let $\xi^\circ_+$ be a component of the annihilator of $\xi$ in
$T^*M$ minus the zero section: $\xi^\circ\smallsetminus 0 =\xi^\circ_+
\sqcup (-\xi^\circ_+ )$.  Assume that $M$ is compact and connected and
that the dimension of $G$ is bigger than 2.  If 0 is not in the image
of the contact moment map $\Psi: \xi^\circ_+ \to \fg^*$ then the
fibers of $\Psi$ are connected and the moment cone $C(\Psi) = \Psi (\xi
^\circ_+)
\cup \{0\}$ is a convex rational polyhedral cone.
\end{theorem}
\begin{remark}
A polyhedral set in $\fg^*$ is the intersection of finitely many
closed half-spaces.  A polyhedral set is rational if the annihilators
of codimension one faces are spanned by vectors in the {\bf integral
lattice} $\ell$ of $\fg$, that is, by vectors in the kernel of $\exp
:\fg \to G$.  The whole space $\fg^*$ is trivially a rational
polyhedral cone. Note that a rational polyhedral cone $C$ in $\fg^*$
is of the form $$ C = \bigcap_i \{ v_i \geq 0\} $$ for some finite
collection of vectors $v_1, \ldots, v_r$ in the integral lattice
$\ell$.
\end{remark}
\begin{remark}
For actions of tori of dimension less than or equal than 2, the fibers
of the corresponding moment maps need not be connected.  For action of
two-dimensional tori the moment cone need not be convex.  In fact, it
is easy to construct an example of an effective 2-torus action on an
overtwisted 3-sphere so that the image cone is not convex.  It is also
easy to construct examples of moment maps for actions of 2-tori and
circles with non-connected fibers (the convexity result for circles is
trivial).  See \cite{IJL}.
\end{remark}

Theorem~\ref{main theorem} extends known convexity results for
Hamiltonian torus actions on symplectic manifolds.  Such results have
a long history.  Atiyah \cite{at:convexity} and, independently,
Guillemin and Sternberg \cite{gu-st:convexity} proved that for
Hamiltonian torus actions on compact symplectic manifolds the image of
the moment map is a rational polyhedron and that the fibers of the
moment map are connected.  The assumption of compactness of the
manifold has been subsequently weakened by de Moraes and Tomei
\cite{MT}, by Prato
\cite{pr:convexity}, by Hilgert, Neeb, and Plank
\cite{hi-ne-pl:convexitycoadjoint} using the methods of
\cite{co-da-mo:geometrie} and by Lerman, Meinrenken, Tolman and Woodward 
\cite{LMTW} to the point where it is enough to assume that the
moment map is proper as a map from a symplectic manifold $M$ to a
\underline{convex} open subset $U$ of the dual of the Lie algebra $\fg^*$.  
The conclusion is that the fibers of the moment map are connected and
that the intersection of the image of the moment map with $U$ a convex
locally polyhedral set.  Note that the hypotheses of Theorem~\ref{main
theorem} only guarantee that the moment map $\Psi : \xi _+^\circ \to
\fg^*$ is proper as a map into $\fg^* \smallsetminus \{0\}$, which is
certainly not convex.

Theorem~\ref{main theorem} is  a direct generalization of a convexity
theorem of Banyaga and Molino \cite{BM}:
\begin{theorem}[Banyaga - Molino]\label{thm_BM}
Let $(M, \xi)$ be a co-oriented contact manifold with an effective
contact action of a torus $G$ preserving the co-orientation.  Assume
that $M$ is compact and connected, that the dimension of $G$ is
bigger than 2 and that $\dim M + 1 = 2 \dim G$. Then the moment cone
$C(\Psi)$ is a convex rational polyhedral cone.
\end{theorem}

\begin{remark}
It is easy to show the hypotheses of the Banyaga - Molino theorem guarantee
that the image of the moment map does not contain the origin:
\begin{lemma}
Let $(M, \xi)$ be a co-oriented contact manifold with an effective
action of a torus $G$ preserving the contact structure and its
co-orientation. Let $\alpha$ be an invariant contact form with $\ker
\alpha = \xi$ and let $\Psi_\alpha :M \to \fg^*$ be the corresponding
moment map.  If $\dim M + 1 = 2 \dim G$ then $\Psi_\alpha (x) \not =
0$ for any $x\in M$.
\end{lemma}

\begin{proof}
Suppose not. Then for some point $x\in M$ the orbit $G\cdot x$ is
tangent to the contact distribution.  Therefore the tangent space
$\zeta _x := T_x (G\cdot x)$ is isotropic in the symplectic vector
space $(\xi _x , \omega_x)$ where $\omega_x = d\alpha _x|_\xi$.

We now argue that this forces the action of $G$ not to be effective.
More precisely we argue that the slice representation of the connected
component of identity $H$ of the isotropy group of the point $x$ is
not effective. The group $H$ acts on $\xi_x$ preserving the symplectic
form $\omega_x$ and preserving $\zeta_x= T_x (G\cdot x)$.  Since
$\zeta_x$ is isotropic, $\xi_x = (\zeta_x^\omega/ \zeta_x) \oplus
(\zeta_x \times \zeta_x ^*)$ as a symplectic representation of $H$.
Here $\zeta_x ^\omega$ denotes the symplectic perpendicular to
$\zeta_x$ in $(\xi_x, \omega_x)$.  Note that since $G$ is a torus, the
action of $H$ on $ \zeta_x$ is trivial.  Hence it is trivial on $\zeta_x^*$.

Observe next that the dimension of the symplectic vector space $V =:
\zeta_x^\omega/ \zeta_x$ is $\dim \xi_x - 2 \dim \zeta_x = 
\dim M -1 - 2 (\dim G - \dim H) = (\dim M - 1) - (\dim M + 1) + 2 \dim H
 = 2 \dim H - 2$.  
On the other hand, since $H$ is a compact connected Abelian group acting
symplecticly on $V$, its image in the group of symplectic linear
transformations $\SP (V)$ lies in a maximal torus $T$ of a maximal
compact subgroup of $\SP (V)$.  The dimension of $T$ is $\dim V/2 =
\dim H -1$.  Therefore the representation of $H$ on $V$ is not
faithful.  Since the fiber at $x$ of the normal bundle of $G\cdot x$
in $M$ is $(T_xM/\xi _x) \oplus (\xi_x/\zeta_x) \simeq \R \oplus (V \oplus
\zeta_x^*)$, the slice representation of $H$ is not faithful.
Consequently the action of $G$ in not effective in a neighborhood of an
orbit $G\cdot x$.  Contradiction.
\end{proof}

\end{remark}

\begin{remark}
As far as we know the paper \cite{BM} is not published.  It is a
revision of \cite{BM2}, which is not widely available.
Theorem~\ref{thm_BM} is cited without proof in \cite{B}.  Providing an
independent and easily accessible proof of Theorem~\ref{thm_BM} is
one of the  motivations for this paper.
\end{remark}

\begin{remark}
We do not know if the condition that no orbit is tangent to the
contact distribution is necessary for Theorem~\ref{main theorem} to
hold.
\end{remark}
\subsection*{A note on notation}  

Throughout the paper the Lie algebra of a Lie group denoted by a
capital Roman letter will be denoted by the same small letter in the
fraktur font: thus $\fg$ denotes the Lie algebra of a Lie group $G$
etc.  The vector space dual to $\fg$ is denoted by $\fg^*$. The identity element of a Lie group is denoted by 1.  The natural 
pairing between $\fg$ and $\fg^*$ will be denoted by 
$\langle \cdot, \cdot \rangle$.

When a Lie group $G$ acts on a manifold $M$ we denote the action by an
element $g\in G$ on a point $x\in G$ by $g\cdot x$; $G\cdot x$ denotes the $G$-orbit of $x$ and so on.  The
vector field induced on $M$ by an element $X$ of the Lie algebra $\fg$
of $G$ is denoted by $X_M$.  The isotropy group of a point $x\in M$ is
denoted by $G_x$; the Lie algebra of $G_x$ is denoted by $\fg_x$ and
is referred to as the isotropy Lie algebra of $x$.  We recall that
$\fg_x = \{ X \in \fg\mid X_M (x) = 0\}$.

If $P$ is a principal $G$-bundle then $[p, m]$ denotes the point in the
associated bundle $P\times _G M = (P\times M)/G$ which is the orbit of
$(p,m) \in P\times M$.

\subsection*{Acknowledgments} I thank Stephanie Alexander for commenting on a
 draft of this manuscript.

\section{Torus actions on contact manifolds}

We now proceed with a proof of Theorem~\ref{main theorem}.  The
methods we use is a mixture of the ideas from
\cite{co-da-mo:geometrie} and \cite{LMTW}.

Recall that $M$ denotes a compact connected manifold with an effective
action of a torus $G$ ($\dim G >2$) preserving a co-oriented contact
distribution $\xi$.  Choose a $G$-invariant contact form $\alpha$ with
$\ker \alpha = \xi$.  Let $\Psi _\alpha : M \to \fg^*$ be the
corresponding moment map; it is defined by equation (\ref{eq1}).
Recall also that we assume that $0\not \in \Psi_\alpha (M)$.  Note that this
condition amounts to saying that no orbit of $G$ is tangent to the
contact distribution $\xi $; thus it is a condition on a contact
distribution and not on a particular choice of a contact form
representing the distribution.

Next fix an inner product on the dual of the Lie algebra
$\fg^*$. Since $\Psi _\alpha (x) \not = 0$ for all $x$ we can define a
new contact form $\alpha'$ by 
$$
\alpha' _x := \frac{1}{||\Psi_\alpha (x)||} \alpha _x .
$$ 
Then the corresponding moment map $\Psi_{\alpha'}$ satisfies
$||\Psi_{\alpha'} (x)|| = 1$ for all $x\in M$.  We assume from now on
that we have chosen an invariant contact form $\alpha$ in such a way
that the corresponding moment map $\Psi _\alpha$ sends $M$ to the unit
sphere $S:= \{ f \in \fg^* \mid ||f || = 1\}$.

\begin{lemma} \label{lemma1a}
Let $(M, \xi)$ be a co-oriented contact manifold with an effective
contact action of a torus $G$. Assume that no orbit of $G$ is tangent
to the contact distribution. Let $\alpha $ be a $G$-invariant contact
form defining $\xi$ normalized so that the image of $M$ under the
corresponding moment map $\Psi_\alpha $ lies in the unit sphere $S$ in
$\fg^*$.  Let $H \subset \fg^*$ be an open half-space, i.e., suppose
that for some $0\not =v\in \fg$ we have $H =\{ f \in \fg^* \mid
\langle f, v\rangle > 0\}$.

For any connected component $N$ of  $\Psi_\alpha \inv (H)$, the fibers of $\Psi _\alpha |_N$ are connected.
\end{lemma} 

\begin{lemma} \label{lemma1b}
Let $M$, $\xi$, $G$, $\alpha$ and $\Psi_\alpha$ be as in
Lemma~\ref{lemma1a} above.  Let $H$ be an open half-space and $N$ a
component of $\Psi_\alpha \inv (H)$.  Then $\Psi_\alpha (N)$ is a
convex rational polyhedral subset of $H \cap S \subset \fg^*$ with
open interior.
\end{lemma}

\begin{remark}
A subset $W$ of the unit sphere $S = \{ f\in \fg^* \mid ||f || = 1\}$
 is convex iff there is a convex cone $C\subset
\fg^*$ (with the vertex at the origin) so that $W = S \cap C$.
Equivalently, $W$ is convex if for any two points $x,y\in W$ there is
a geodesic of length $\leq \pi$ connecting $x$ to $y$ and lying
entirely in $W$.

A subset $W$ of $S$ (respectively of $H \cap S$) is rational polyhedral if 
there exist vectors $v_1, \ldots v_k$ in the integral lattice 
$\ell = \ker \{ \exp :\fg \to G\}$ such that
$$
W = \{f \in S \mid \langle f, v_i\rangle \geq 0, \quad 1\leq i \leq k\}
$$ 
(respectively if $W = \{f \in S \cap H  \mid \langle f, v_i\rangle \geq 0, \quad 1\leq i \leq k\}$).
\end{remark}

\begin{proof}[Proof of Lemmas~\ref{lemma1a} and \ref{lemma1b}]
Consider the symplectization $(M\times \R, d(e^t \alpha))$ of $(M, \alpha)$.  As usual $t$ denotes the coordinate on $\R$.  The contact action of $G$ on $M$ extends trivially to a Hamiltonian action on the symplectization.  The corresponding moment map $\Phi :M\times \R \to \fg^*$ is given by 
$$
\Phi (x, t) = e^t \Psi_\alpha (x). 
$$ 
The symplectic manifold $(N\times \R, d(e^t \alpha)|_{N\times \R})$
is a symplectization of $(N, \alpha |_N)$. The manifold $N\times \R$
is a connected symplectic manifold with a Hamiltonian action of $G$,
the map $\Phi_N := \Phi|_{N\times \R}$ is a corresponding moment map
for the action of $G$.  Moreover, it has the following two properties:
\begin{enumerate}
\item $\Phi_N (N\times \R)$ is contained in the convex open subset $H$ 
of $\fg^*$;
\item $\Phi_N :N\times \R \to H$ is proper.  
\end{enumerate}
Therefore Theorem~4.3 of \cite{LMTW} applies.  We conclude that the
fibers of $\Phi_N$ are connected and that the image $\Phi_N (N\times
\R)$ is convex.

Next, since the action of the torus $G$ on $M$ is effective, it is
free on a dense open subset of $M$.  This is a consequence of the
principal orbit type theorem and the fact that $G$ is abelian.
Consequently the action of $G$ on $N\times \R$ is free on a dense open
subset.  Hence the image $\Phi_N (N\times \R)$ has non-empty interior.
Also, since $M$ is compact and $G$ is abelian, the number of subgroups
of $G$ that occur as isotropy groups of points of $M$ is finite.
Therefore not only does [LMTW, Theorem~4.3] imply that $\Phi_N (N\times \R)$
is the intersection a locally polyhedral subset of $\fg^*$ with the
open half-space $H$, but that in fact $\Phi_N (N\times \R) = \Phi
(N\times \R)$ is a polyhedral cone.
\end{proof}

\begin{lemma}\label{lemma2}
Let $M$, $G$, $\alpha$ and $\Psi_\alpha$ be as in Lemma~\ref{lemma1a}
above.  Define an equivalence relation $\sim$ on $M$ by declaring the
equivalence classes to be the connected components of the fibers of
the moment map $\Psi _\alpha $.  Let $\bar{M} = M/\sim$.

Then $\bar{M}$ is a compact path connected space and the moment map
$\Psi_\alpha : M\to \fg^*$ descends to a continuous map $\bar{\Psi} :
\bar{M} \to S$, where as before $S$ is the unit sphere in $\fg^*$
centered at 0.

Moreover, $\bar{M}$ is a length space and $\bar{\Psi}: \bar{M} \to S$
is  locally an isometric embedding.  More precisely, for any open
half-space $H$ and any connected component $N$ of $\bar{\Psi}\inv (H)$
the map $\bar{\Psi}|_N : N \to S$ is an isometric embedding.
\end{lemma}

Our proof of Lemma~\ref{lemma2} uses length spaces, the notion that is
due to Gromov \cite{Gromov1, Gromov2}.  We therefore briefly summarize
the relevant facts.  Our treatment follows a book by D. Burago,
Yu. Burago and S. Ivanov \cite{BBI}.

\subsection{Digression: length structures and length spaces}

Let $X$ be a topological space.  Consider a class $\A$ of continuous
paths in $X$ which is closed under restrictions, concatenations and
reparameterizations.  Suppose that there is a map $L: \A \to [0,
\infty]$ (``the length'') satisfying the following conditions for any
curve $\gamma :[a, b] \to X$ in $\A$:
\begin{itemize}
\item[a] $L(\gamma) = L(\gamma|_{[a, c]}) + L(\gamma |_{[c, b]})$ for any $c \in (a, b)$.
\item[b] The function $L_t :=  L(\gamma|_{[a, t]})$ is a continuous function of $t\in [a, b]$.

\item[c] If $\varphi: [c, d] \to [a, b]$ is monotone and continuous, then $L (\gamma) = L (\gamma \circ \varphi)$. 

\item[d] If a sequence of curves $\gamma_i \in \A$ converges to $\gamma$ uniformly, then $L(\gamma) \leq \liminf L(\gamma_i)$.

\item[e] If $U\subset X$ is a proper open subset, and $p\in U$ is a point 
then the number 
$$
\inf \{L(\gamma) \mid \gamma: [a, b] \to X, \gamma \in \A, \gamma (a) = p, 
\gamma (b) \not \in U\}
$$
is positive.
\end{itemize}
\begin{definition}
The triple $(X, \A, L)$ where $X$ is a topological space, $\A$ is a
class of continuous curves in $X$ and $L:\A \to [0, \infty]$ is a map
satisfying the conditions above is called a {\bf length structure}.
\end{definition}

Suppose $(X, \A, L)$  is a length structure.   Suppose  that for any two points $x,y \in X$ there is a path $\gamma \in \A$ starting at $x$ and ending  at $y$. We then define the {\bf distance} $d_L :X \times X \to [0, \infty]$ by 
$$
d_L (x, y) = \inf \{ L(\gamma) \mid \gamma :[a, b] \to X, \gamma (a) = x, \gamma (y) = b, \gamma \in \A\}.
$$
One can check that if $d_L (x,y) <\infty$ for all $x,y\in X$ then $d_L$ is  a metric.

Suppose $(X, d)$ is a metric space. Then we can take $\A$ to be the set of rectifiable paths and $L =L_d:\A \to [0, \infty]$ to be the length functional.  Then $(X, \A, L)$ is a length structure.  Note that in general $d_L (x,y) \geq d(x,y)$ for $x,y\in X$.  If $d_L = d$ then $(X, d)$ is called a {\bf length space}.  A unit sphere $S$ in a normed finite dimensional vector space with the standard metric induced by the embedding is an example of a length space.

\begin{definition}
Let $(X, \A, L)$ be a length structure.   Let $\gamma: [a, b] \to X$ be a curve in $\A$.  It is a {\bf geodesic} if for any $c, d\in [a,b]$ with $|c -d|$ sufficiently small $L(\gamma|_{[c, d]}) = d_L (\gamma(c), \gamma (d))$. 
\begin{remark} We think of geodesics as maps, not as subsets.
Also, from now on all geodesics are parameterized by arc length.
\end{remark}
 
\end{definition} 
If $(X, d)$ is a compact connected metric space then a version of
Hopf-Rinow theorem holds, and so any two points of $X$ can be
connected by a geodesic. See, for example, Proposition~3.7 in
\cite{BH}.  This ends our digression on length spaces.

\begin{proof}[Proof of Lemma~\ref{lemma2}]
It is clear that $\bar{M}$ is a compact path-connected topological
space and that the moment map $\Psi_\alpha: M\to \fg^*$ descends to a
continuous map $\bar{\Psi}: \bar{M} \to S = \{ ||f || = 1\} $  Moreover, by
Lemmas~\ref{lemma1a} and \ref{lemma1b}, for any open half-space
$H\subset \fg^*$ and any component $Z$ of $\bar{\Psi}\inv (H)$, the
map $\bar{\Psi}: Z \to S\cap H$ is a topological embedding which is a
homeomorphism on an open dense set.

This gives us a way to define a length structure on $\bar{M}$: We
define the class $\A$ to be the set of all curves $\bgamma: [a, b] \to
\bar{M}$ such that $\bar{\Psi} \circ \bgamma$ is a rectifiable curve
in the unit sphere $S$.  For $\bgamma \in \A$ we set $L(\bgamma) = L_S
(\bar{\Psi} \circ \bgamma)$ where $L_S$ is the length functional on
the rectifiable curves in the sphere defined by the standard metric.
Let $d_L$ be the corresponding metric on $\bar{M}$.  Then, since for
any half-space $H$ and any component $Z$ of $\bar{\Psi}\inv (H)$ the
set $\bar{\Psi} (Z)$ is convex in the sphere $S$, the map $ \bar{\Psi}
:Z \to S$ is an {\em isometric} embedding.  Thus $\bar{\Psi}: \bar{M}
\to S$ is locally an isometric embedding.
\end{proof}

\begin{corollary}\label{cor1} Let $\bar{M}$, $\bar{\Psi}$ and $S$ be as in Lemma~\ref{lemma2}
If $\bgamma$ is a geodesic in $\bar{M}$ then $\bar{\Psi} \circ \bgamma$ is a geodesic in $S$.
\end{corollary}
\begin{remark}
Since $\bPsi$ is a local isometry it maps geodesics in $\bM$ to
geodesics in the unit sphere $S$ of the {\bf same length}.  In
particular if the end points of a (nonconstant) geodesic $\bgamma$ in
$\bM$ are sent by $\bPsi$ to the same point in the sphere, then $\bPsi
\circ \bgamma$ multiply covers a great circle and consequently the
length of $\bgamma$ is an integer multiple of $2\pi$.
\end{remark}

We emphasize that Lemmas~\ref{lemma1a} and \ref{lemma1b} can be
restated for the induced map $\bar{\Psi} : \bar{M} \to S$ of
Lemma~\ref{lemma2} as follows:

\begin{lemma}\label{lemma1a'} 
For any open half-space $H$ and any connected component $N$ of 
$\bar{\Psi}\inv (H)$ the map $\bar{\Psi}|_N \to S$ is an isometric embedding.
\end{lemma}

\begin{lemma}\label{lemma1b'}  
 For any open half-space $H$ and any connected component $N$ of
 $\bar{\Psi}\inv (H)$ the set $\bar{\Psi} (N)$ is a convex polyhedral
 subset of the sphere $S$ with non-empty interior.
\end{lemma}

A consequence of Lemmas~\ref{lemma1a'} and \ref{lemma1b'} we get:

\begin{corollary} \label{cor'}
Let $\bar{\Psi}: \bar{M} \to S$ be as in Lemma~\ref{lemma2}. Suppose the points $x_1, x_2 \in \bar{M}$ lie in the same connected component of $\bar{\Psi}\inv (H)$  for some  open half-space $H$.  

If $\bPsi (x_1) = \bPsi (x_2)$ then $x_1 = x_2$.  If $\bPsi (x_1) \not
= \bPsi (x_2)$ then there is a geodesic $\bgamma$ in $\bM$ connecting
$x_1$ to $x_2$.  Moreover we may choose $\bgamma$ such that $\bPsi
\circ \bgamma$ is a geodesic in $S$ lying entirely in the half-space
$H$ and connecting $\bPsi (x_1)$ and $\bPsi (x_2)$.
\end{corollary}
As a consequence of Lemma~\ref{lemma2} we get:

\begin{corollary} \label{cor2}
Any two points in $\bM$ can be connected by a short geodesic, i.e.,
for any two points $x, y\in \bM$ there is a geodesic $\bgamma$ with
$\bgamma (0) = x$ and $\bgamma (d) = y$ where $d$ is the distance
between $x$ and $y$ (recall that all geodesics are parameterized by
arc length).

\end{corollary}

\begin{remark}
Such a geodesic in $\bM$ need not be unique.  For example consider the
unit co-sphere bundle $M$ in the cotangent bundle of a flat torus $G$.
Then $M = G\times S$, $\Psi :G\times S \to S \subset \fg^*$ is the
projection and $\bM$ is the unit sphere $S$.  In this case for any
point $x\in \bM = S$ there are infinitely many geodesics of length
$\pi$ connecting $x$ and $-x$.
\end{remark}

The following lemma  uses the notation above.
\begin{lemma}\label{lemmaS}
Suppose $x_1, x_2$ are two points in $\bM$ connected by a path
$\bgamma$ with the property that $\bPsi \circ \bgamma$ lies entirely
in some open half-space $H$.  Then the points $x_1$, $x_2$ lie in the
same connected component of $\bPsi \inv (H)$.
\end{lemma}
\begin{proof}
The image of $\bgamma$ lies in a connected component of $\bPsi \inv (H)$.
\end{proof}
Lemma~\ref{lemma3} below is the main technical tool for proving the
connectedness of fibers of moment maps.

\begin{lemma}\label{lemma3}
Let $\bPsi :\bM \to S$ be as in Lemma~\ref{lemma2}.  Suppose
$\bgamma_1$, $\bgamma_2$ are two distinct geodesics in $\bM$ with
$\bgamma_1 (0) = \bgamma_2 (0)$, and suppose that $\bPsi \circ \gamma_1$
and $\bPsi\circ \gamma_2$ trace out two {\em distinct} great circles in the
unit sphere $S$.  Then $\bgamma_2(0) = \bgamma_2 (2\pi)$ (and so
$\bgamma_1 (0) = \bgamma _1 (2\pi)$).
\end{lemma}
\begin{remark}
Note that the assumption  $\dim G>2$ is crucial for the lemma to make sense.
\end{remark}

\begin{proof}[Proof of Lemma~\ref{lemma3}]
The idea of the proof is to show that there is an open half-space $H$
containing $\bPsi (\bgamma_2 (0))$ such that $\bgamma _2 (0)$ and
$\bgamma _2 (2\pi)$ lie in the same connected component of $\bPsi \inv
(H)$. For then by Corollary~\ref{cor'} $\bgamma_2(0) = \bgamma_2
(2\pi)$.

Given a path $\bgamma _i$ in $\bM$ we write $\gamma_i$ for the path
$\bPsi \circ \bgamma_i$ in $S$.

Since by assumption the geodesics $\gamma_1$ and $\gamma _2$ trace out
two distance great circles in $S$, $\gamma _1 (\frac{\pi}{2}) \not =
\pm \gamma _2(\frac{\pi}{2})$.  On the other hand we clearly have
$\gamma _1 (0) = -\gamma _1 (\pi) = -\gamma _2 (\pi)$, $\gamma _1
(2\pi) = \gamma_2 (2\pi) = \gamma_1 (0)$, $\gamma_1 (\frac{3\pi}{2}) =
- \gamma _1 (\frac{\pi}{2})$ and $\gamma_2 (\frac{3\pi}{2}) = - \gamma
_2 (\frac{\pi}{2})$.

Since $\gamma _1 (\frac{\pi}{2}) \not = \pm \gamma _2(\frac{\pi}{2})$,
there is an open half-space $H_1$ containing the points $\gamma _1
(0)$, $\gamma _1 (\frac{\pi}{2})$ and $\gamma _2 (\frac{\pi}{2})$.  By
Lemma~\ref{lemmaS}, $\bgamma_1 (\frac{\pi}{2})$ and $\bgamma_2
(\frac{\pi}{2})$ lie in the same connected component of $\bPsi \inv
(H_1)$ as $\bgamma _1(0)$.  By Corollary~\ref{cor'} there a geodesic
$\bsigma_1$ in $\bM$ connecting $\bgamma_1 (\frac{\pi}{2})$ to
$\bgamma_2 (\frac{\pi}{2})$ such that $\sigma_1 := \bPsi \circ
\bsigma_1$ traces out a short geodesic connecting $\gamma_1
(\frac{\pi}{2})$ to $\gamma _2 (\frac{\pi}{2})$.

Choose an open half-space $H_2$ containing the points $\gamma_1
(\frac{\pi}{2})$, $\gamma_2 (\frac{\pi}{2})$ and $\gamma_1 (\pi)
=\gamma_2 (\pi)$. Note that by construction $\bgamma_1
(\frac{\pi}{2})$ is connected to $\bgamma _2 (\frac{\pi}{2})$ by
$\bsigma _1$, $\bgamma_1 (\frac{\pi}{2})$ is connected to $\bgamma_1
(\pi)$ by a piece of $\bgamma_1$ and $\bgamma_2 (\frac{\pi}{2})$ is
connected to $\bgamma_2 (\pi)$ by a piece of $\bgamma_2$.  By
Lemma~\ref{lemmaS} $\bgamma _1 (\pi)$ and $\bgamma _2(\pi)$ lie in the
same connected component of $\bPsi \inv (H_2)$.  By
Corollary~\ref{cor'} we have $\bgamma_1 (\pi)= \bgamma _2 (\pi)$.

Choose a half-space $H_3$ containing $\gamma_1 (\pi)$, $\gamma _1
(\frac{\pi}{2})$ and $\gamma _2 (\frac{3\pi}{2})$.  Since $\bgamma _1
(\pi) = \bgamma_2 (\pi)$, since $\bgamma _1 (\pi)$ is connected to
$\bgamma_1 (\frac{\pi}{2})$ by a piece of $\bgamma_1$ and since
$\bgamma_2 (\pi)$ is connected to $\bgamma _2 (\frac{3\pi}{2})$ by a
piece of $\bgamma_2$, $\bgamma_1 (\frac{\pi}{2})$ and $\bgamma_2
(\frac{3\pi}{2})$ lie in the same connected component of $\bPsi \inv
(H_3)$.  By Corollary~\ref{cor'} there a geodesic $\bsigma_2$ in $\bM$
connecting $\bgamma_1 (\frac{\pi}{2})$ to $\bgamma_2 (\frac{3\pi}{2})$
such that $\sigma_2 := \bPsi \circ \bsigma_2$ traces out a short
geodesic connecting $\gamma_1 (\frac{\pi}{2})$ to $\gamma _2
(\frac{3\pi}{2})$.

Finally choose a half-space $H_4$ containing $\gamma _1(0) = \gamma _2
(0) = \gamma _2 (2\pi)$, $\gamma _1 (\frac{\pi}{2})$ and $\gamma _2
(\frac{3\pi}{2})$.  Arguing as above we see that $\bgamma _2 (0)$ and
$\bgamma _2 (2\pi)$ lie in the same connected component of $\bPsi \inv
(H_4)$.  Hence, by Corollary~\ref{cor'}, $\bgamma_2 (0) = \bgamma _2
(2\pi)$.
\end{proof}

\begin{lemma}\label{lemma4}
 The fibers of the orbit moment map $\bPsi :\bM \to S$ are connected,
 i.e., $\bPsi$ is an embedding.
\end{lemma}

\begin{proof}
Suppose $x_1, x_2 \in \bM$ are two points with $\bPsi (x_1) = \bPsi
(x_2)$.  We want to show that $x_1 = x_2$.  Suppose not.  Then the
distance $d$ between $x_1$ and $x_2$ is positive.  Let $\bgamma_1$ be
a short geodesic connecting $x_1$ and $x_2$, so that $\bgamma _1 (0) =
x_1$ and $\bgamma _1 (d) = x_2$.  Then $\gamma _1 := \bPsi \circ
\bgamma _1$ is a geodesic in the unit sphere $S$ starting and ending
at $\gamma_1 (0)$.  Therefore $\gamma _1$ multiply covers a great
circle in $S$ (and so $d$ is an integer multiple of $2\pi$).
 
Suppose that we can construct a geodesic $\bgamma_2$ connecting $x_1$
to $x_2$ so that $\gamma_2 := \bPsi \circ \bgamma _2$ covers a great
circle distinct from the one covered by $\gamma_1$.  Then by
Lemma~\ref{lemma3} $\bgamma _1 (0) = \bgamma _1 (2\pi)$ contradicting
the choice of $\bgamma_1$ as a short geodesic.

Now we construct $\bgamma_2$ with the required properties.  Pick an
open half-space $H$ containing $\gamma_1 (0)$.  Let $N$ denote the
connected component of $\bPsi \inv (H)$ containing $x_1$. By
Lemma~\ref{lemma1b'} the set $\bPsi (N)$ is convex with nonempty
interior.  Pick a point $y$ in $N$ so that $\bPsi (y)$ is not in the
image of the geodesic $\gamma _1$. By Corollary~\ref{cor'} there is a
geodesic $\bsigma$ connecting $x_1$ to $y$ with the image of $\sigma
:= \bPsi \circ \bsigma$ lying entirely in $H$.  Let $\btau$ be a short
geodesic connecting $y$ to $x_2$.  If the image of $\tau := \bPsi
\circ \btau$ lies entirely in a half-space containing $\bPsi (x_2)$
and $\bPsi (y)$ then by Lemma~\ref{lemmaS} we have $x_1 = x_2$.

Otherwise $\tau$ traces out a long geodesic connecting $\bPsi (y)$ to $\bPsi (x_2) =\gamma _1 (0)$.  If $\btau$ passes through $x_1$ then the piece of $\btau$ starting at $x_1$ and ending at $x_2$ is the desired geodesic $\bgamma_2$.  If $\btau$ does not pass through $x_1$, concatenate $\bsigma$ with $\btau$.  The concatenation $\bgamma_2$ is the desired geodesic. 
\end{proof}
\begin{lemma}\label{lemma_conv}
The image of the orbit moment map $\bPsi :\bM \to S$ is convex.
\end{lemma}
\begin{proof}
Suppose $f_1, f_2$ are two points in the image of $\bPsi$.  Then either $f_1$ and $f_2$ lie in some open half-space $H$ or $f_1 = - f_2$.  In the former case, by Lemma~\ref{lemma4}, $N= \bPsi\inv (H)$ is connected.  Hence, by Lemma~\ref{lemma1b'}, $\bPsi (N) = H \cap \bPsi (\bM )$ is convex and consequently $\bPsi (\bM )$ is convex.  

In the latter case we argue as follows.  The sets $\bPsi \inv (f_i)$, $i=1,2$ consists of single points; denote these points by $x_i$.  Connect $x_1$ and $x_2$ by a short geodesic $\bgamma$.  Then the image of $\gamma = \bPsi \circ \bgamma$ contains an arc of a great circle in $S$ passing through $f_1$ and $f_2 = -f_1$ (in fact it follows from the proof of Lemma~\ref{lemma3} that the image of $\gamma$ is exactly such an arc). 
\end{proof}

\begin{lemma}\label{lemma_image_cone}
Let $\Psi _\alpha : M \to \fg^*$ be a moment map as in
Lemma~\ref{lemma1a}.  The corresponding moment cone $C(\Psi)$ is a
rational convex polyhedral cone.  That is either $C(\Psi) = \fg^*$ or
there exist vectors $v_1, \ldots, v_k$ in the integral lattice $\ell $
of the torus $G$ such that $$ C(\Psi) = \bigcap _i \{v_i \geq 0\}.  $$
\end{lemma}

\begin{proof}
By Lemmas~\ref{lemma1b'} and \ref{lemma4} for any open half-space $H$ of $\fg^*$ there exist vectors $v_1, \ldots, v_r$ in the integral lattice $\ell$ ($r$ depends on $H$) such that
$$
C(\Psi) \cap H = \left( \bigcap _i \{v_i \geq 0\} \right) \cap H.
$$
Moreover, we may and will assume that the set of $v_i$'s is minimal.
Thus no $v_i$ is strictly positive on $C(\Psi) \cap H$.  Since the
moment cone is a cone on a compact set, there exist finitely many
half-spaces $H^1, \ldots, H^s$ such that $\bigcup_\beta H^\beta$
contains $C(\Psi)\smallsetminus \{0\}$.  For each such half-space
$H^\beta$, let $v_1^\beta, \ldots, v_{r(\beta)}^\beta$ be the minimal
set of integral vectors so that
$$
C(\Psi) \cap H^\beta = 
	\left( \bigcap _i \{v_i ^\beta \geq 0\} \right) \cap H^\beta.
$$
We claim that
$$
C(\Psi) = \bigcap _{i, \beta} \{v_i^\beta \geq 0\} .
$$
As a first step we argue that for any $i, \beta$ we have 
$$
C(\Psi) \subset \{v_i^\beta \geq 0\} .
$$
By choice of $v_i ^\beta$ there exists a point $x\in C(\Psi) \cap
H^\beta$ such that $v_i^\beta (x) = 0$ (since $x\in H^\beta$, $x\not =
0$).  Suppose there exists a point $y\in C(\Psi)$ with $v_i^\beta (y)
<0$.  Since $C(\Psi)$ is convex, $tx + (1-t)y \in C(\Psi)$ for all
$t\in [0, 1]$. On the other hand $v_i ^\beta (tx + (1-t)y ) = (1-t)
v_i^\beta (y) < 0 $ for all $t\in [0, 1)$.  Since $H^\beta$ is open
there is $\epsilon >0$ so that $tx + (1-t)y \in H^\beta $ for all
$t\in (\epsilon, 1]$.  Therefore for all $t\in (\epsilon, 1)$ we have
$$ 
tx + (1-t)y \in H^\beta \cap C(\Psi) \subset \{v_i^\beta \geq 0\},
$$ 
a contradiction.  We conclude that 
$$ 
C(\Psi) \subset \bigcap _{i,\beta} \{v_i^\beta \geq 0\}.  
$$
%Note that since $C(\Psi) \subset \bigcup H^\beta \cup \{0\}$, $C(\Psi) %\subset \left(\bigcup H^\beta \cup \{0\}\right) \cap \bigcap _{i, \beta} %\{v_i^\beta \geq 0\}$.

Next we argue that the reverse inclusion holds as well: $\bigcap _{i,
\beta} \{v_i^\beta \geq 0\} \subset C(\Psi)$.  By construction for
each $\beta$ $$ C(\Psi) \cap H^\beta = \left( \bigcap _i \{v_i ^\beta
\geq 0\} \right) \cap H^\beta.  $$ Since $\bigcup _\beta H^\beta\cup
\{0\}$ covers the image cone $C(\Psi)$, we have
\begin{eqnarray*}
C(\Psi)= C(\Psi)\cap (\bigcup _\beta H^\beta\cup \{0\}) &=& \{0\} \cup \bigcup _\beta (C(\Psi)\cap H^\beta )\\
&=& \bigcup _\beta \left( (\bigcap _i \{ v_i^\beta \geq 0\} \cap (H^\beta\cup \{0\}\right)\\
&\supseteq& \left ( \bigcap_{i, \beta}  \{ v_i^\beta \geq 0\} \right) \cap \left(\bigcup _\beta H^\beta\cup \{0\}\right) 
\end{eqnarray*}
Therefore
\begin{equation}\label{eq*}
C(\Psi)= \left ( \bigcap_{i, \beta}  \{ v_i^\beta \geq 0\} \right) \cap (\bigcup _\beta H^\beta\cup \{0\}) .
\end{equation}
Finally, since $\bigcap_{i, \beta}  \{ v_i^\beta \geq 0\} $ is closed and convex, its  intersection with the unit sphere $S\cap \bigcap_{i, \beta}  \{ v_i^\beta \geq 0\} $ is closed and connected. On the other hand 
\begin{equation}\label{eq**}
S \cap \bigcap_{i, \beta}  \{ v_i^\beta \geq 0\}  = \left ( S \cap \bigcap_{i, \beta}  \{ v_i^\beta \geq 0\}  \cap (\bigcup _\beta H^\beta) \right) \sqcup S\cap  \left( \bigcap_{i, \beta}  \{ v_i^\beta \geq 0\}  \smallsetminus (\bigcup _\beta H^\beta)\right).
\end{equation}
It follows from (\ref{eq*}) and (\ref{eq**}) that the set $S \cap \bigcap_{i, \beta}  \{ v_i^\beta \geq 0\} $ is a disjoint union of two closed sets.   Therefore the set $S\cap \left(\bigcap_{i, \beta}  \{ v_i^\beta \geq 0\}  \smallsetminus \cup _\beta H^\beta \right)$ is empty.  We conclude that
$$
C(\Psi) = \bigcap_{i, \beta}  \{ v_i^\beta \geq 0\}  \cap (\bigcup _\beta H^\beta\cup \{0\})= \bigcap_{i, \beta}  \{ v_i^\beta \geq 0\}.
$$
\end{proof}

\end{document}